\documentclass{article}
\usepackage{amsmath,amssymb,amsthm}
\newtheorem{Thm}{Theorem}
\newtheorem{Lem}{Lemma}
\newtheorem{Prop}{Proposition}
\theoremstyle{definition}

\begin{document}

\title{Uniqueness of positive solutions to semilinear elliptic equations with double power nonlinearities}
\author{Shinji Kawano}
\date{}
\maketitle

\begin{abstract}

We consider the problem
\begin{equation}
 \begin{cases}
   \triangle u+f(u)=0  & \text{in  $\mathbb{R}^n$},\\
   \displaystyle \lim_{\lvert x \rvert \to \infty} u(x)  =0, \label{Intro2}
 \end{cases}  
\end{equation}
where  
\begin{equation*}
f(u)=-\omega u + u^p - u^q, \qquad  \omega>0, ~~q>p>1. 
\end{equation*}
It is known that a positive solution to \eqref{Intro2} exists if and only if $F(u):=\int_0^u f(s)ds >0$ for some $u>0$.
Moreover, Ouyang and Shi in 1998 found that the solution is unique if $f$ satifies furthermore the condition
that $\tilde{f}(u):= (uf'(u))'f(u)-uf'(u)^2 <0$ for any $u>0$.
In the present paper we remark that this additional condition is unnecessary.
\end{abstract}

\section{Introduction}

We shall consider a boundary value problem
\begin{equation}
 \begin{cases}
   u_{rr}+ \dfrac{n-1}{r}u_r+f(u)=0 & \text{for $r>0$}, \\
   u_r(0)=0, \\
   \displaystyle \lim_{r \to \infty} u(r) =0, \label{b}
 \end{cases} 
\end{equation} 
where $n \in \mathbb{N}$ and
\begin{equation*}
f(u)=-\omega u + u^p - u^q, \qquad \omega>0, ~~q>p>1. 
\end{equation*}              
The above problem arises in the study of 
\begin{equation}
 \begin{cases}
   \triangle u+f(u) =0  & \text{in  $\mathbb{R}^n$},\\
   \displaystyle \lim_{\lvert x \rvert \to \infty} u(x)  =0. \label{a}
 \end{cases}  
\end{equation}
Indeed, the classical work of Gidas, Ni and Nirenberg~\cite{G1,G2} tells us that any positive solution to \eqref{a} is radially symmetric. 
On the other hand, for a solution $u(r)$ of \eqref{b}, $v(x):=u(\lvert x \rvert)$ is a solution to \eqref{a}.  

The condition to assure the existence of positive solutions to \eqref{a} (and so \eqref{b}) was given by 
Berestycki and Lions~\cite{B1} and Berestycki, Lions and Peletier~\cite{B2}:

\begin{Prop}
A positive solution to \eqref{b} exists if and only if 
\begin{equation}
F(u):=\int_0^u f(s)ds >0, \qquad \text{for some} \quad u>0.       \label{existence}
\end{equation} 
\end{Prop}     

Uniqueness of positive solutions to \eqref{b} had long remained unknown.
Finally in 1998 Ouyang and Shi~\cite{OS} proved uniqueness for \eqref{b} with $f$ satisfying the additional condition
(See also Kwong and Zhang~\cite{KZ}):

\begin{Prop}
If $f$ satisfies furthermore the following condition, then the positive solution is unique;
\begin{equation}
\tilde{f}(u):= (uf'(u))'f(u)-uf'(u)^2 <0, \qquad \text{for any} \quad u>0.      \label{unique}
\end{equation}     
\end{Prop}

Following is the main result of the present paper:
\begin{Thm}
If the nonlinearity $f$ satisfies the existence condition~\eqref{existence}, then the uniqueness condition~\eqref{unique}
is automatically fulfilled.
\end{Thm}   \label{thm}

This paper is organized as follows. 
In section~2 we give a straightforward proof of the theorem.
In section~3 we give an alternative proof of the theorem, in which an interesting technical lemma is used.
In section~4 we explain the technical lemma.   

\section{Proof of Theorem~\ref{thm}.}

\begin{Lem}
The existence condition~\eqref{existence} is equivalent to
\begin{equation*}
\omega < \omega_{p,q},
\end{equation*}
 where 
\begin{equation*}
\omega_{p,q}=\dfrac{2(q-p)}{(p+1)(q-1)} \left[ \dfrac{(p-1)(q+1)}{(p+1)(q-1)} \right] ^{\frac{p-1}{q-p}}. 
\end{equation*}
(See Ouyang and Shi~\cite{OS} and the appendix of Fukuizumi~\cite{Fukuizumi}.) 
\end{Lem}

\begin{Lem}
The uniqueness condition~\eqref{unique} is equivalent to
\begin{equation*}
\omega < \eta_{p,q},
\end{equation*}
where
\begin{equation*}
\eta_{p,q}=\dfrac{q-p}{q-1}\left[ \dfrac{p-1}{q-1}\right]^{\frac{p-1}{q-p}}.
\end{equation*}
\end{Lem}

The proofs of these Lemmas are nothing but straightforward calculation and shall be omited.

\begin{proof}[Proof of Theorem~\ref{thm}.]
It is apparent that
\begin{equation*}
0<\omega_{p,q}<\eta_{p,q},
\end{equation*}
which asserts the theorem.
\end{proof}

\section{Alternative Proof of Theorem~\ref{thm}.}

The following is the key lemma. The proof of this lemma is given in the next section as a corollary of more general statement.

\begin{Lem}
The existence condition~\eqref{existence} is equivalent to the following condition;
\begin{equation}
\tilde{F}(u)= (uf(u))'F(u)-uf(u)^2 <0, \qquad \text{for any} \quad u>0.      \label{exitilde}
\end{equation} 

The uniqueness condition~\eqref{unique} is equivalent to the following condition;
\begin{equation}
f(u) >0, \qquad \text{for some} \quad u>0.       \label{uniquetilde}
\end{equation} 
\end{Lem}    \label{lemtilde}

This Lemma means that the relation between $f$ and $\tilde{f}$ is parallel to the relation between $F$ and $\tilde{F}$.

\begin{proof}[Alternative Proof of Theorem~\ref{thm}.]
From Lemma~\ref{lemtilde} it is enough to show that if $F$ has positive parts (i.e. the existence condition~\eqref{existence}) 
then $f$ has positive parts (i.e. the condition~\eqref{uniquetilde}). 
The contraposition of this statement is clear by the monotonicity of the integral.
\end{proof}

\section{Classification of double power nonlinear functions}

To state the main result of this section, from now on we mean a more general function
\begin{equation}
 f(u)=-au^p+bu^q-cu^r,\qquad \text{for}  ~u>0,    \label{general}
\end{equation} 
where $a,b,c>0$ and $p<q<r$, by the same notation~$f$.
We also steal the former notation 
\begin{equation}
\tilde{f}(u)= (uf'(u))'f(u)-uf'(u)^2 .      \label{tilde}
\end{equation}    
Following is the statement.
\begin{Thm}
There only can occur the following three cases;
\begin{itemize}
\item[(a)]~~~~~ $a<b\dfrac{r-q}{r-p}\left[ \dfrac{b(q-p)}{c(r-p)} \right]^{\frac{q-p}{r-q}}$ $\Longleftrightarrow$ $f$ has positive parts  $\Longleftrightarrow$ $\tilde{f}$ remains negative.
\item[(b)]~~~~~ $a=b\dfrac{r-q}{r-p}\left[ \dfrac{b(q-p)}{c(r-p)} \right]^{\frac{q-p}{r-q}}$ $\Longleftrightarrow$ $f$ has just one zero $\Longleftrightarrow$ $\tilde{f}$ has just one zero.
\item[(c)]~~~~~ $a>b\dfrac{r-q}{r-p}\left[ \dfrac{b(q-p)}{c(r-p)} \right]^{\frac{q-p}{r-q}}$ $\Longleftrightarrow$ $f$ remains negative $\Longleftrightarrow$ $\tilde{f}$ has positive parts.
\end{itemize}
\end{Thm}    \label{th}

\begin{proof}
The statement with respect to $f$ is trivial.
We obtain from the definition~\eqref{tilde} that
\begin{equation*}
\tilde{f} = -ab(q-p)^2 u^{q+p-1} + ca(r-p)^2 u^{r+p-1} - bc(r-q)^2 u^{r+q-1} .
\end{equation*}
This is in the form of~\eqref{general} and use the result with respect to $f$.  
\end{proof}

At last we give the proof of Lemma~\ref{lemtilde}, which concludes this paper.

\begin{proof}[Proof of Lemma~\ref{lemtilde}.]
For the first part, we consider
\begin{equation*}
F(u)= -\frac{\omega}{2}u^2+\frac{u^{p+1}}{p+1}-\frac{u^{q+1}}{q+1}, \qquad  \omega>0, ~~q>p>1,
\end{equation*}
which is in the form of~\eqref{general} and use the result (a) of Theorem~\ref{th}. 

For the second part, we consider
\begin{equation*}
f(u)=-\omega u + u^p - u^q, \qquad  \omega>0, ~~q>p>1,
\end{equation*}
which is in the form of~\eqref{general} and use the result (a) of Theorem~\ref{th}. 
\end{proof}


\begin{thebibliography}{99}
\bibitem{B1}
  H. Berestycki and P. L. Lions, 
  Nonlinear scalar field equation, I.,
  \textit{Arch. Rat. Meth. Anal.} \textbf{82} (1983), 313-345.
\bibitem{B2}
  H. Berestycki, P. L. Lions and L. A. Peletier, 
  An ODE approach to the existence of positive solutions for semilinear problems in $\mathbb{R}^N$,
  \textit{Indiana University Math. J.} \textbf{30} (1981), 141-157.
\bibitem{Fukuizumi}
  R. Fukuizumi,
  Stability and instability of standing waves for nonlinear Schr\"{o}dinger equations,
  \textit{Tohoku Mathematical Publications}, \textbf{No.25}, 2003.
\bibitem{G1}
  B. Gidas, W. M. Ni and L. Nirenberg,
  Symmetry and related properties via the maximal principle,
  \textit{Comm. Math. Phys.} \textbf{68} (1979), 209-243. 
\bibitem{G2}
  B. Gidas, W. M. Ni and L. Nirenberg,
  Symmetry of positive solutions of nonlinear elliptic equations in $\mathbb{R}^n$, Mathematical analysis and applications, Part A,
  \textit{Adv. in Math. Suppl. stud.} \textbf{7a}, Academic press, 1981, 369-402.
\bibitem{Kawano}
  S. Kawano,
  A remark on the uniqueness of positive solutions to semilinear elliptic equations with double power nonlinearities, 
  preprint.  
\bibitem{KZ}
  M. K. Kwong and L. Zhang, 
  Uniqueness of positive solutions of $\triangle u+f(u)=0$ in an annulus,
  \textit{Differential Integral Equations} \textbf{4} (1991), 583-599.
\bibitem{M}
  T. Mizumachi,
  Uniqueness of positive solutions to a scalar field equation with a double power nonlinearity,
  preprint.
\bibitem{OS}
  T. Ouyang and J. Shi,
  Exact multiplicity of positive solutions for a class of semilinear problems,
  \textit{J. Differential Equations} \textbf{146} (1998), 121-156.
\bibitem{PS}  
  L. A. Peletier and J. Serrin,
  Uniqueness of positive solutions of semilinear equations in $\mathbb{R}^n$,
  \textit{Arch. Rat. Mech. Anal.} \textbf{81} (1983), 181-197.
\bibitem{WW}
  J. Wei and M. Winter,
  On a cubic-quintic Ginzburg-Landau equation with global coupling,
  \textit{Proc. Amer. Math. Soc.} \textbf{133} (2005), 1787-1796.
\end{thebibliography}
\end{document}